\documentclass[11pt]{article}

\usepackage{graphicx}
\usepackage{amsmath}
\usepackage{amsthm}
\usepackage{amssymb}
\usepackage[mathscr]{eucal}

\title{Every mapping class group is generated by $6$ involutions}
\author{Tara E. Brendle and
Benson Farb \thanks{The first author is partially supported by a
VIGRE postdoc under NSF grant number 9983660 to Cornell University.
The second author is supported in part by NSF grants DMS-9704640
and DMS-0244542.}}

\theoremstyle{plain}
\newtheorem{theorem}{Theorem}

\newtheorem{lemma}[theorem]{Lemma}
\newtheorem{corollary}[theorem]{Corollary}

\newtheorem{question}[theorem]{Question}

\def\proof{{\bf {\medskip}{\noindent}Proof. }}

\def\endproof{$\diamond$ \bigskip}

\def\title{\em}

\def\bar{\overline}

\newcommand\R{\mbox{\bf R}}
\newcommand\Z{\mbox{\bf Z}}

\newcommand\mgb{{\rm Mod_{g,b}}}
\newcommand\mgo{{\rm Mod_{g,1}}}
\newcommand\mgz{{\rm Mod_{g,0}}}
\newcommand\Mod{{\rm Mod}}
\newcommand\surf{S_{g,b}}
\newcommand\sgb{S_{g,b}}
\newcommand\sgz{S_{g,0}}

\DeclareMathOperator\SL{SL}

\DeclareMathOperator\Sp{Sp}

\begin{document}
\maketitle
\begin{abstract}
Let $\mgb$ denote the mapping class group of a surface of genus
$g$ with $b$ punctures.  Luo asked in \cite{Lu} if there is a
universal upper bound, independent of genus, for the number of
torsion elements needed to generate $\mgb$.  We answer Luo's
question by proving that $3$ torsion elements suffice to generate
$\mgz$.  We also prove the more delicate result that
there is an upper bound, independent of genus, not only for the
number of torsion elements needed to generate $\mgb$ but also for
the order of those elements.  In particular, our main result is
that $6$ involutions (i.e. orientation-preserving diffeomorphisms
of order two) suffice to generate $\mgb$ for every genus $g \geq
3, b = 0$ and $g \geq 4, b = 1$.
\end{abstract}

\section{Introduction}
Let $S_{g,b}$ denote a closed, oriented surface of genus $g$ with
$b$ punctures, and let $\mgb$ denote its mapping class
group, which is the group of homotopy classes of
orientation-preserving homeomorphisms preserving the set of punctures.
We shall frequently abuse
terminology by confusing an individual homeomorphism with its
mapping class in $\mgb$.

We begin with a brief survey of some known generating sets for
$\mgb$ possessing various properties.  Dehn \cite{De} produced a
finite set of generators of $\mgz$, proving that $2g(g-1)$ Dehn
twists suffice for $g \geq 3$. Lickorish \cite{Li} improved on
this result by giving a generating set for $\mgz$ consisting of
$3g - 1$ twists for any $g \geq 1$. Humphries \cite{Hu} then
showed that a certain subset of Lickorish's set consisting of $2g
+ 1$ twists suffices to generate $\mgz$ and that this is in fact
the minimal number of twist generators for $\mgz$ (here again $g
\geq 1$). Johnson \cite{Jo} later proved that Humphries' set also
generates $\mgo$.

If one allows generators other than twists, smaller generating
sets can be obtained. Lickorish \cite{Li} noted that $\mgz$ can be
generated with $4$ elements, 3 of which are twists and 1 of which
has finite order.  N. Lu \cite{Lu} found a
generating set with $3$ elements, $2$ of which have finite order.
Wajnryb \cite{Wa} proved that for $g \geq 1$, $b =
0,1$, the group $\mgb$ can be generated by $2$ elements, one of
which has finite order.

The problem of finding small generating sets, torsion generating sets,
and generating sets of {\em involutions} (i.e.~elements of order two) 
is a classical one, and has been studied extensively, especially for
finite groups (see, e.g., \cite{DT} for a survey).  In 1971, Maclachlan
\cite{Mac} proved that $\mgz$ is generated by torsion elements, and
deduced from this that moduli space ${\cal M}_g$ is simply-connected as a
topological space.  These results were
later extended to $\mgb, g\geq 3,b\geq 1$ by Patterson \cite{Pa}.  The
question of generating mapping class groups by involutions\footnote{We
remind the reader that the only involutions under consideration are
orientation-preserving.} was first investigated by McCarthy and
Papadopoulos \cite{MP}. Among other results, they proved that for $g
\geq 3$, $\mgz$ is generated by infinitely many conjugates of a certain
involution.

Luo \cite{Luo}, using work of Harer \cite{Ha}, described the first
finite set of involutions which generate $\mgb$ for $g \geq 3, b
\geq 0$. The order of his generating set depends on both $g$ and
$b$; in particular, his set consists of $12g + 2$ involutions when
$b = 0,1$.  Luo also gives torsion generators, not necessarily
involutions, for all other cases except $g = 2, b = 5k + 2$.  It
should be noted that the existence of a generating set for $\mgz$
consisting of $4g + 4$ torsion elements follows directly from a
lemma of Birman \cite{Bi} (see Lemma~\ref{jb} below).  In his
paper (\cite{Luo}, \S1.4), Luo poses the question of whether there
is a universal upper bound, independent of $g$ and $b$, for the
number of torsion elements needed to generate $\mgb$.  Our first
result is a positive answer to Luo's question for $b=0$.

\begin{theorem}[Three torsion elements generate]
\label{torsion}
For each $g \geq 1$, the group $\mgz$ is generated by $3$ elements of
finite order.
\end{theorem}

As we will see in Section~\ref{thmtwo}, at least one of the three
torsion generators we give has order depending on $g$.  Subsequent to
the original posting of this paper, M. Korkmaz \cite{Ko1} has shown
that $\mgb, b=0,1$ is generated by two elements, each of order $4g+2$.

Finding a set of generators whose orders are {\em universally} bounded
and whose cardinality is also universally bounded is more delicate,
especially if one wants a generating set consisting of
involutions. Our main theorem addresses this.

\begin{theorem}[Six involutions generate]
\label{main} For $g \geq 3, b = 0$, and for $g \geq 4, b = 1$, the
group $\mgb$ is generated by $6$ involutions.
\end{theorem}

In \cite{Ka}, Kassabov builds on our method to 
extend and improve this result to the case 
$b>1$.  Further, in some cases (e.g. $g\geq 8$) he proves that $4$
involutions suffice to generate $\mgb$.

In \S\ref{section:final}, we note that Theorem~\ref{main} implies that 
$\mgb$ is the quotient of a $6$-generator Coxeter group.

\bigskip
\noindent
{\bf Remarks. }
\begin{enumerate}
\item Since $\Mod_{1,0}=\Z/4\Z \ast_{\displaystyle \Z/2\Z} \Z/6\Z$ and
${\rm H}_1(\Mod_{2,0},\Z)=\Z/10\Z$, it is easy to see that
Theorem \ref{main} does not extend to the genus $g=1$ or $g=2$ cases.

\item As $\mgz$ surjects onto the integral symplectic group
$\Sp(2g,\Z)$, it follows from Theorems~\ref{torsion}
and~\ref{main} that for all $g\geq 3$, the group $\Sp(2g,\Z)$ is
generated by $3$ torsion elements, and also by $6$ involutions.

\item If one allows {\em orientation-reversing} involutions, it is
possible to use Theorem \ref{torsion}, together with other arguments, to
prove that the {\em extended mapping class group}, which includes
orientation-reversing mapping classes, is generated by $5$ involutions.

\end{enumerate}

It would be interesting to obtain the exact bounds for Theorem 
\ref{main}.  It is easy to see the lower bound of $3$; we do not know
how to improve on this bound.

It is a pleasure to thank Dieter Kotschick, Mustafa Korkmaz, and Dan
Margalit for their valuable suggestions.  We also thank Ian Agol for
suggesting that Remark $3$ above might be possible, Nathan Broaddus for pointing out the connection with Coxeter groups, and Martin
Kassabov who pointed out a redundancy in our original $7$-involution 
generating set, thus reducing the conclusion in Theorem~\ref{main} 
to $6$ involutions.

\section{Proof that $3$ torsion elements generate}
\label{thmtwo}

In this section we shall only consider the case of the closed surface,
i.e., $b = 0$.

We begin with a lemma of Birman \cite{Bi},
and include an adapted version of her proof for
completeness.  For a simple closed curve $x$ in $\surf$, let
$T_x$ denote the (right-handed) Dehn twist about $x$.

\begin{lemma}[Two torsion elements for a twist] \label{jb}
Let $x$ be a simple closed curve which is nonseparating in $\sgz,
\hskip .05in g\geq 1$.
Then $T_x$ can be written as the product of two torsion elements.
\end{lemma}

\begin{figure}
$$
\setlength{\unitlength}{0.05in}
\begin{picture}(0,0)(0,11)
\put(15.5,45.5){ {\bf {\footnotesize $\alpha_1$} } }
\put(9,38.5){ {\bf {\footnotesize $\beta_1$} } }
\put(19,30){ {\bf {\footnotesize $\gamma_1$} } }
\put(32.5,38.5){ {\bf {\footnotesize $\alpha_2$} } }
\put(22.5,38){ {\bf {\footnotesize $\beta_2$} } }
\put(23.5,25.5){ {\bf {\footnotesize $\gamma_2$} } }
\put(39.5,21.5){ {\bf {\footnotesize $\alpha_3$} } }
\put(32.5,27.7){ {\bf {\footnotesize $\beta_3$} } }
\put(24,19){ {\bf {\footnotesize $\gamma_3$} } } \put(-2,40){ {\bf
{\footnotesize $\alpha_g$} } } \put(-1.5,29.5){ {\bf
{\footnotesize $\beta_g$} } } \put(8.5,27){ {\bf {\footnotesize
$\gamma_{g-1}$} } } \put(28,45){ {\bf {\footnotesize $\rho_1$}} }
\end{picture}
\includegraphics[width=2in]{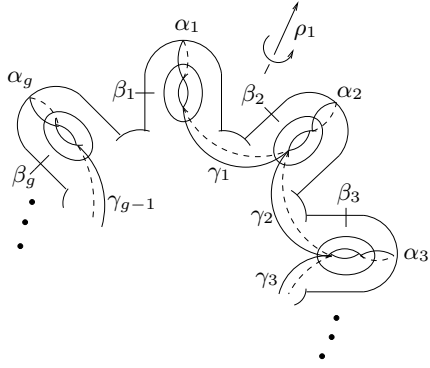}
$$
\caption{A circular embedding of the surface together with
Lickorish's $3g-1$ twist generators and an involution $\rho_1$
which is rotation by $\pi$ about the axis indicated.} \label{Li2}
\end{figure}

\proof  We will find it convenient to fix a ``circular embedding''
of the surface $\sgb$, as seen in Figure~\ref{Li2}.  Referring to
the curves of the same figure, we define
\begin{equation}
\label{eq:definingQ}
Q = T_{\alpha_g}
T_{\beta_g} (T_{\gamma_{g-1}} T_{\beta_{g-1}}) (T_{\gamma_{g-2}}
T_{\beta_{g-2}}) \cdots (T_{\gamma_1} T_{\beta_1}) T_{\alpha_1}
\end{equation}
and
\begin{equation}
\label{eq:definingS}
S = Q T_{\alpha_1}^{-1} = T_{\alpha_g} T_{\beta_g} (T_{\gamma_{g-1}}
T_{\beta_{g-1}}) (T_{\gamma_{g-2}} T_{\beta_{g-2}}) \cdots (T_{\gamma_1}
T_{\beta_1})
\end{equation}

We say that an ordered set ${c_1, \ldots, c_n}$ of simple closed curves on
$\surf$  forms an {\em n-chain} if the geometric
intersection $(c_k, c_{k+1}) = 1$ for $1 \leq k \leq n-1$ and $(c_k,
c_l) = 0$ if $|k - l| \geq 2$.  If $n$ is odd, the boundary of a regular
neighborhood of any $n$-chain has two components $d_1$ and $d_2$; if $n$
is even the boundary has one component $d$.  The so-called {\em chain
relation} in $\mgb$ tells us that for a given $n$-chain ${c_1, \ldots,
c_n}$, if $n$ is odd we have
$$(T_{c_1} T_{c_2} \cdots T_{c_n})^{n + 1}
= T_{d_1} T_{d_2}$$ and if $n$ is even, we have
$$(T_{c_1}T_{c_2} \cdots T_{c_n})^{2n + 2} = T_d$$

The curves defining the twists of $Q$ and $S$ form a $(2g + 1)$-chain
and a $(2g)$-chain, respectively.  However, any boundary curve of a
regular neighborhood of either chain is null-homotopic in the closed
surface.  Hence the chain relation tells us that $Q$ and $S$ have orders
which divide $2g + 2$ and $4g + 2$, respectively, in $\mgz$.  In fact,
Birman notes in her original proof that an ``ugly but routine''
calculation establishes that these values are precisely the orders of
the two elements, with explicit calculations recorded in \cite{BH} (see
also \cite{HK} for a shorter proof).

Given any nonseparating simple closed curve $x$, there is a
homeomorphism $h$ such that $h (\alpha_1) = x$.
Recall that for $h \in \mgb$ and a simple closed curve $c$ contained in
$\surf$, we have
\begin{equation}
\label{eq:conj}
h T_c h^{-1} =T_{h(c)}
\end{equation}
Then we have
\begin{eqnarray*}
T_x &=& h T_{\alpha_1} h^{-1} \\
&=& h (S^{-1} Q) h^{-1} \\
&=& (h S^{-1} h^{-1}) (h Q h^{-1})
\end{eqnarray*}
which gives the desired result.
\endproof


Now to complete the proof of Theorem \ref{torsion}.  It is well
known that $\Mod_{1,0}\cong \SL(2,\Z)$ is generated by an element
of order $4$ and an element of order $6$, so we can assume $g\geq
2$.  Let $\alpha_i$ be as above, let $Q$ and $S$ be defined as in
(\ref{eq:definingQ}) and (\ref{eq:definingS}), and let
$$U=T_{\alpha_1} T_{\alpha_2}^{-1}$$

Wajnryb \cite{Wa} showed that $\mgz, \hskip .05in  g\geq 1$ is
generated by the two maps $U$ and $S$.  In the proof of
Lemma~\ref{jb}, we saw that $S$ has finite order in $\mgz$; hence
it remains to deal with $U$, which is the product of two Dehn
twists.  As above, we have

$$T_{\alpha_1} = S^{-1}
Q$$

From the proof of Lemma~\ref{jb}, we know that $$ T_{\alpha_2} =
h (S^{-1} Q) h^{-1} $$ where $h$ is any map taking the curve
$\alpha_1$ to $\alpha_2$.

Let $\rho_1$ denote the involution which is rotation by $\pi$
about the axis indicated in Figure~\ref{Li2}.  Clearly, $\rho_1
(\alpha_1) = \alpha_2$ , so we can write
\begin{eqnarray*}
U &=& T_{\alpha_1} T_{\alpha_2}^{-1} \\
&=& [S^{-1} Q] [\rho_1 (Q^{-1} S) \rho_1^{-1}]
\end{eqnarray*}
Thus we can generate Wajnryb's two generators, and hence all of
$\mgz$, with three torsion elements: $Q, S$ and $\rho_1$.
\endproof

{\bf Remarks.}
\begin{enumerate}
\item
We can replace the involution $\rho_1$ in the above proof with a
rotation of order $g$.
\item
Wajnryb \cite{Wa} also shows that the two maps $U$ and $S$
generate $\mgo$.  However the proof of Theorem~\ref{torsion} does
not go through for $b = 1$ since Lemma~\ref{jb} fails in this
case.
\end{enumerate}

We also note that, using Wajnryb's two generators $S, U$, we can
generate $\mgb$ when $b = 0,1$ with a set consisting of two
involutions and one element of order $4g + 2$.  Since $U =
T_{\alpha_1} T_{\alpha_2}^{-1}$, and $\rho_1 (\alpha_2) =
\alpha_1$, we have that $U = T_{\alpha_1} (\rho_1
T_{\alpha_1}^{-1} \rho_1)$.  Shifting the parentheses, we have
that $U$ is the product of the involutions $\rho_1$ and its
conjugate $T_{\alpha_1} \rho_1 T_{\alpha_1}^{-1}$.  We will see
the usefulness of this ``parentheses shifting" technique in the
next section.

\section{A universal bound on involutions}

For the remainder of the paper, we assume that $ g \geq 3$ and $b =
0,1$.  For simplicity of exposition, we provide explicit
arguments only for $b=0$.  In the case $b=1$ the arguments are the same,
although some involutions must be replaced with certain conjugates which
move the puncture to a fixed point of the involution.

\subsection{Writing Dehn twists as a product of involutions}
\label{twists}

We begin by recalling the so-called {\em lantern relation} in
$\mgb$. This relation was discovered by Dehn \cite{De} in the
1930s and was rediscovered by Johnson \cite{Jo} over forty years
later.  For convenience, we will use the notation $X$ in addition
to $T_x$ to denote the (right-handed) Dehn twist about the simple
closed curve $x$ in $\surf$.

\begin{figure}
$$
\setlength{\unitlength}{0.05in}
\begin{picture}(0,0)(0,11)
\put(9,26.5){ {\bf $a_1$} }
\put(18,39){ {\bf $a_2$} }
\put(27,26.5){ {\bf $a_3$} }
\put(18,13){ {\bf $a_4$} }
\put(9,35){ {\bf $x_1$} }
\put(27,35){ {\bf $x_3$} }
\put(18,22){ {\bf $x_2$} }
\end{picture}
\includegraphics[width=2in]{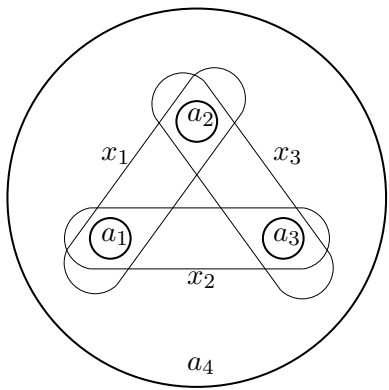}
$$
\caption{The Lantern Relation: $X_1 X_2 X_3 = A_1 A_2 A_3 A_4.$}
\label{lantern}
\end{figure}

Referring to
the four boundary curves $a_1, a_2, a_3$ and $a_4$ of the surface
$S_{0,4}$ together with the interior curves $x_1, x_2$ and $x_3$, as
shown in Figure~\ref{lantern}, the lantern relation is the following
relation amongst the corresponding Dehn twists:
\begin{equation}
\label{eq:lantern}
X_1 X_2 X_3 = A_1 A_2 A_3 A_4
\end{equation}

The symmetry of Figure~\ref{lantern} clearly shows that the twists
on the left-hand side of (\ref{eq:lantern}) can be cyclically
permuted.  Moreover, each twist on the right-hand side commutes
with all other twists in the relation since the corresponding
curves are disjoint from all others in the lantern.  Thus we can
rewrite the relation in the following way, which will be
convenient for our purposes:
\begin{equation}
\label{eq:lantern2}
A_4 = (X_1A_1^{-1}) (X_2 A_2^{-1}) (X_3 A_3^{-1})
\end{equation}

The following lemma improves on a result of Luo \cite{Luo} and
Harer \cite{Ha}; they showed that a Dehn twist about a
nonseparating curve can be written as a product of six
involutions.

\begin{lemma}[Four involutions for a twist]
\label{LH}
Let $c$ be a simple closed curve which is nonseparating in $\surf$, $g
\geq 3$.  Then $T_c$ can be written as the product of four involutions.
\end{lemma}

\proof We begin with the argument of Luo and Harer, which we
include here for completeness.  If $g \geq 3$, we can find a
lantern as in Figure~\ref{lantern} embedded in $\surf$ such that
the given curve $c$ plays the role of $a_4$ and such that the
complement of the lantern is connected and its four boundary
curves are distinct and nonseparating in $\surf$. Let us call such
an embedding a {\em good lantern}. We now observe that for each
$j$, $1 \leq j \leq 3$, the pair $(x_j, a_j)$ consists of two
disjoint, nonseparating simple closed curves such that the
complement of the union $x_j\cup a_j$ is connected in $\surf$ for
each $1\leq j\leq 3$.  Then for each $j$ we can find an involution
$I_j$ such that $I_j (x_j) = a_j$.

By (\ref{eq:conj}) we can write
$$X_j A_j^{-1} = X_j (I_j X_j^{-1} I_j) =
(X_j I_j X_j^{-1}) I_j
$$
\noindent
for each $1\leq j\leq 3$.  Note that each $X_jI_jX_j^{-1}$ is an involution,
so that each $X_j A_j^{-1}, 1\leq j\leq 3$, can be realized as a
product of two conjugate involutions.  Hence
$T_c = A_4$ is the product of six involutions by (\ref{eq:lantern2}).
\begin{figure}
$$
\setlength{\unitlength}{0.05in}
\begin{picture}(0,0)(0,11)
\put(29,29){ {\bf $a_1$} } \put(29,13){ {\bf $a_2$} } \put(29,21){
{\bf $a_3$} } \put(4.5,26){ {\bf $a_4$} } \put(50,26){ {\bf
$a_4$}} \put(11,29.5){ {\bf $x_2$} } \put(11,13){ {\bf $x_1$} }
\put(9.5,21){ {\bf $x_3$} } \put(57,18){ {\bf $J_1$} }

\end{picture}
\includegraphics[width=3in]{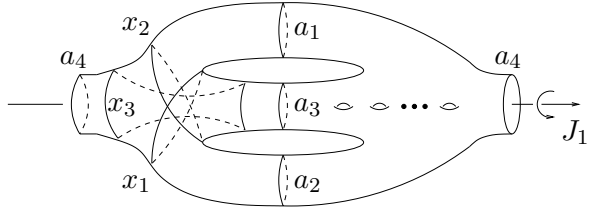}
$$
\caption{The {\em pair swap involution} $J_1$, taking $(a_1,x_1)$ to
$(a_2,x_2)$.  Here the boundary curve
$a_4$ on the left is identified with that on
the right.} \label{twist}
\end{figure}

To improve upon this result, consider the surface shown in
Figure~\ref{twist}.  If the curves labelled $a_4$ are identified
in the obvious way, then the four curves $a_i$ together bound a
good lantern with interior curves $x_j$ as shown.  Given the non-separating curve $c = a_4$, we can find a good lantern with $a_4$ as a boundary component.
We can choose its remaining labels $a_i, x_j, \hskip .05in 1 \leq i,j \leq 3$, in such a way that
there is a homeomorphism of
$S_{g,b}$ to the surface of Figure~\ref{twist} taking this good lantern to the lantern of Figure~\ref{twist}.  It is clear now
that there exists an involution $J_1$ of $S_{g,b}$ which takes the
pair $(a_1, x_1)$ to the pair $(a_2, x_2)$.  We note that the
surface of Figure~\ref{twist} could be embedded in $\R^3$ in such
a way that the involution $J_1$ of the surface is a restriction of
an isometry of $\R^3$; one can simply add a tube connecting the
two ends which encloses the surface of Figure~\ref{twist}. (We
shall abuse notation throughout by using $J_1$ to refer to both
the involution of $\sgb$ as well as its conjugate which is the
involution of the surface of Figure~\ref{twist}.)

In fact, given a particular choice of labels for a good lantern, a
similar process yields $6$ distinct ``pair swaps'' of order $2$,
i.e., involutions taking one pair $(a_i,x_j)$ to another pair
$(a_k, x_l)$, assuming $i \neq k, j \neq l$.  In particular, there
is a conjugate involution $J_2$ taking the pair $(a_1, x_1)$ to
the pair $(a_3, x_3)$.  This is clear if we act on the surface of
Figure~\ref{twists} by the homeomorphism which interchanges the
``tube'' containing the curve $a_2$ with the tube containing $a_3$,
fixing the curve $x_1$ while simultaneously moving the curve $x_3$
to the current position of $x_2$.



We can now rewrite (\ref{eq:lantern2}) as
\begin{equation}
\label{eq:fourinvs}
T_c= A_4 = [ (X_1 I_1 X_1^{-1}) I_1 ] [ J_1 (X_1 I_1 X_1^{-1}) I_1 J_1 ] [
J_2(X_1 I_1 X_1^{-1}) I_1 J_2 ]
\end{equation}
and hence $T_c$ is a product of four involutions.
\endproof

\subsection{Proof of the main result: $6$ involutions generate}
\label{proof}

\begin{figure}
$$
\setlength{\unitlength}{0.05in}
\begin{picture}(0,0)(0,11)
\put(15.5,45){ {\bf $\alpha$} }
\put(13.5,33.5){ {\bf $\beta$} }
\put(23,29){ {\bf $\gamma$} }
\put(23,44){ {\bf $R_g$} }
\end{picture}
\includegraphics[width=2in]{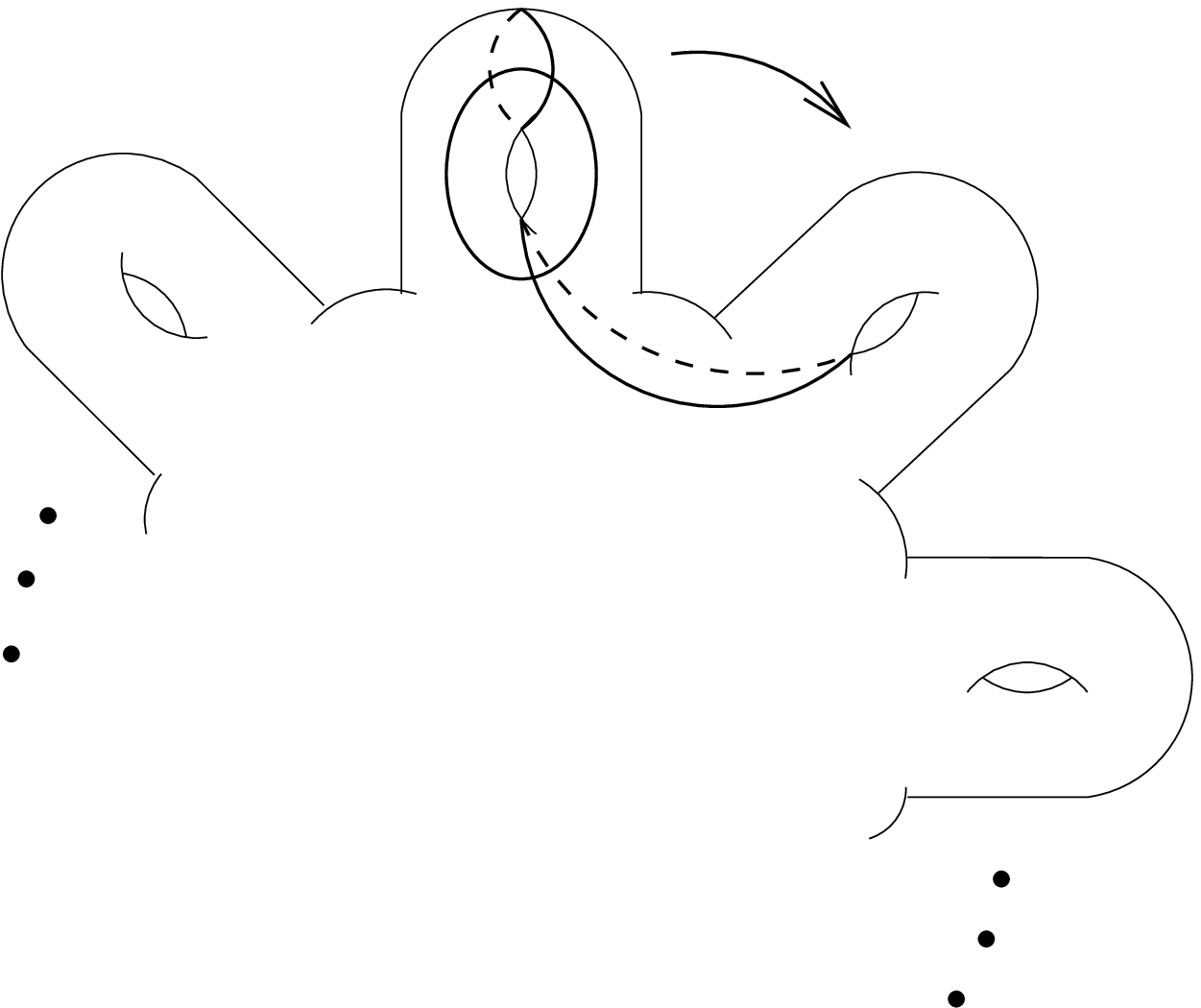}
$$
\caption{A generating set for $\mgb$.}
\label{gens}
\end{figure}

Lickorish \cite{Li} proved that Dehn twists about the $3g -1$
curves given in Figure~\ref{Li2} suffice to generate $\mgb$. In
the same paper, he also described a generating set for $\mgb$
consisting of four elements, one of which is not a Dehn twist.  We
shall now give a rigorous description of such a generating set.
Let $R_g \in \mgb$ denote clockwise rotation of $\surf$ by
$\frac{2 \pi}{g}$ about the axis which is perpendicular to the
plane of the page and intersects the surface twice in the center
of our ``circle of handles''.  Applying Equation (\ref{eq:conj})
to $R_g$ conjugating the twists $T_\alpha, T_\beta$ and $T_\gamma$
(referring to the curves of Figure~\ref{gens}), we see that the
group generated by $R_g$ together with these three twists contains
$$\{R_g^mT_\alpha R_g^{-m}, R_g^mT_\beta R_g^{-m}, R_g^mT_\gamma
R_g^{-m}: 0\leq m \leq g-1\}$$ Hence the four elements $T_\alpha,
T_\beta, T_\gamma, R_g$ generate all of Lickorish's twist generators and
thus they generate all of $\mgb$.  Clearly, one could replace any of the
three twists generators $T_c, \hskip .05in c \in \{ \alpha, \beta,
\gamma \}$, with $R_g^k T_c R_g^{-k} = T_{R_g^k(c)}$ for any $k \in {\bf
Z}$ and the resulting set still generates $\mgb$.

\begin{figure}
\vskip .2in
$$
\setlength{\unitlength}{0.05in}
\begin{picture}(0,0)(0,11)
\put(27,50){ {\bf $L_2$} } \put(15,53.5){ {\bf $L_1$} }
\put(28.5,43){ {\bf $\rho_1$} } \put(19.5,48) { {\bf $\rho_2$} }
\end{picture}
\includegraphics[width=2in]{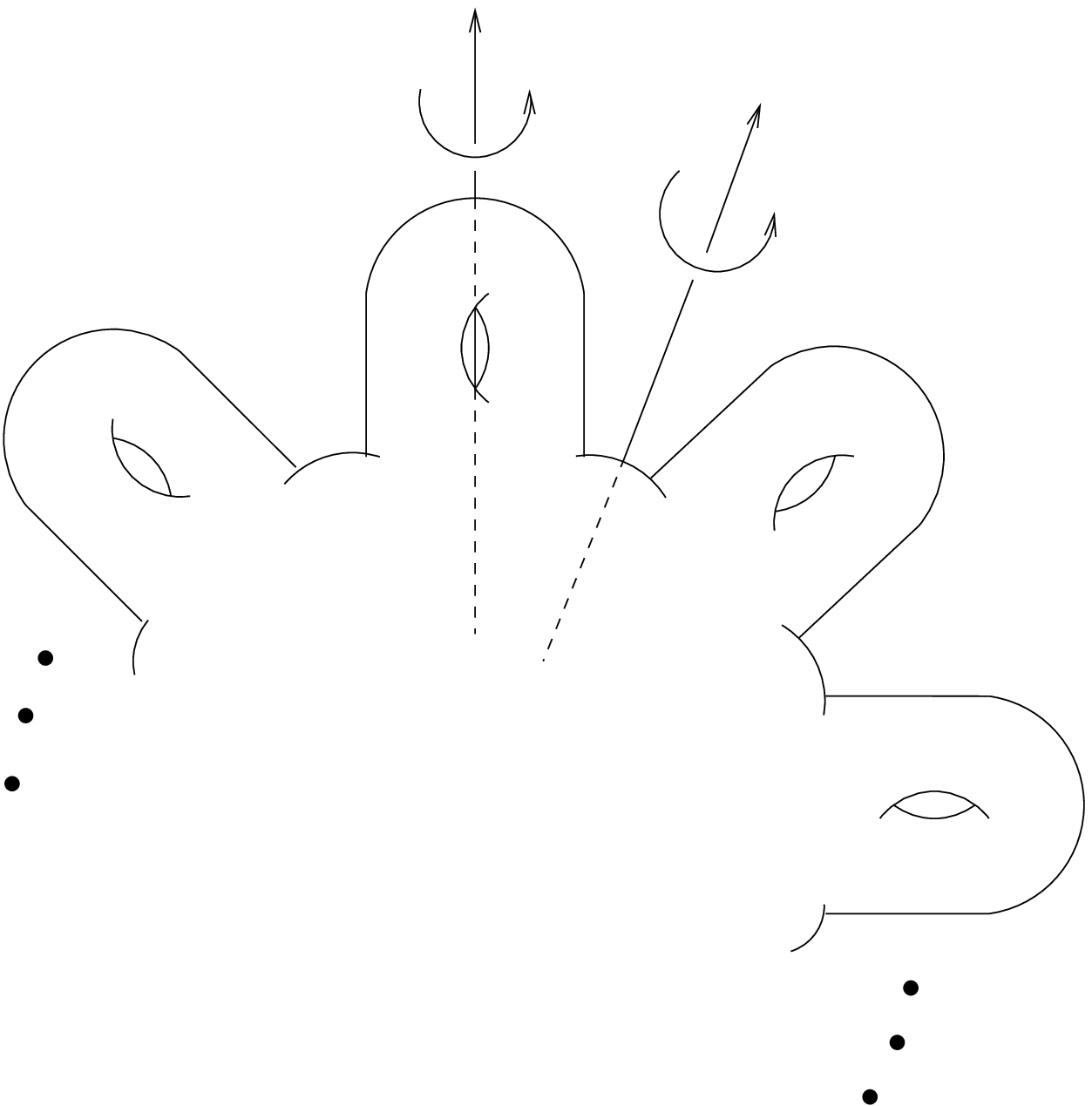}
$$
\caption{Two involutions generating $R_g$.} \label{tr}
\end{figure}

Consider the rotation $R_g$.  For $i = 1,2,$ let $\rho_i$ denote rotation in $\R^3$ by $\pi$ about the axis
labelled $L_i$, as shown in Figure~\ref{tr}.
Note that if the genus $g$ is even, then $L_1$ intersects the
surface in six points and $L_2$ in two points.  If the genus $g$
is odd, then $L_1$ and $L_2$ each intersect the surface in four
points. In both cases, we observe that
\begin{equation}
\label{eq:rotate} R_g = \rho_1 \rho_2
\end{equation}

Combining (\ref{eq:rotate}) with Lemma~\ref{LH}, together with the
fact proved above that $\mgb$ is generated by $R_g$ and the three
twists $T_\alpha, T_\beta, T_\gamma$, we have that $\mgb$ can be
generated by $2+3\cdot 4=14$ involutions.

\begin{figure}
\vskip .2in
$$
\setlength{\unitlength}{0.05in}
\begin{picture}(0,0)(0,11)
\put(12,45.5){ {\bf $\alpha = a_1$} }

\put(13,20){ {\bf $\beta = a_3 $} }

\put(9.5,31.5){ {\bf $\gamma = a_4$} }

\put(33, 39){ {\bf $x_1$} }

\put(27,45.5){ {\bf $\rho_1$} }

\put(24,25){ {\bf $a_2$} }

\put(24,21){\line(1,0){7}}
\end{picture}
\includegraphics[width=2in]{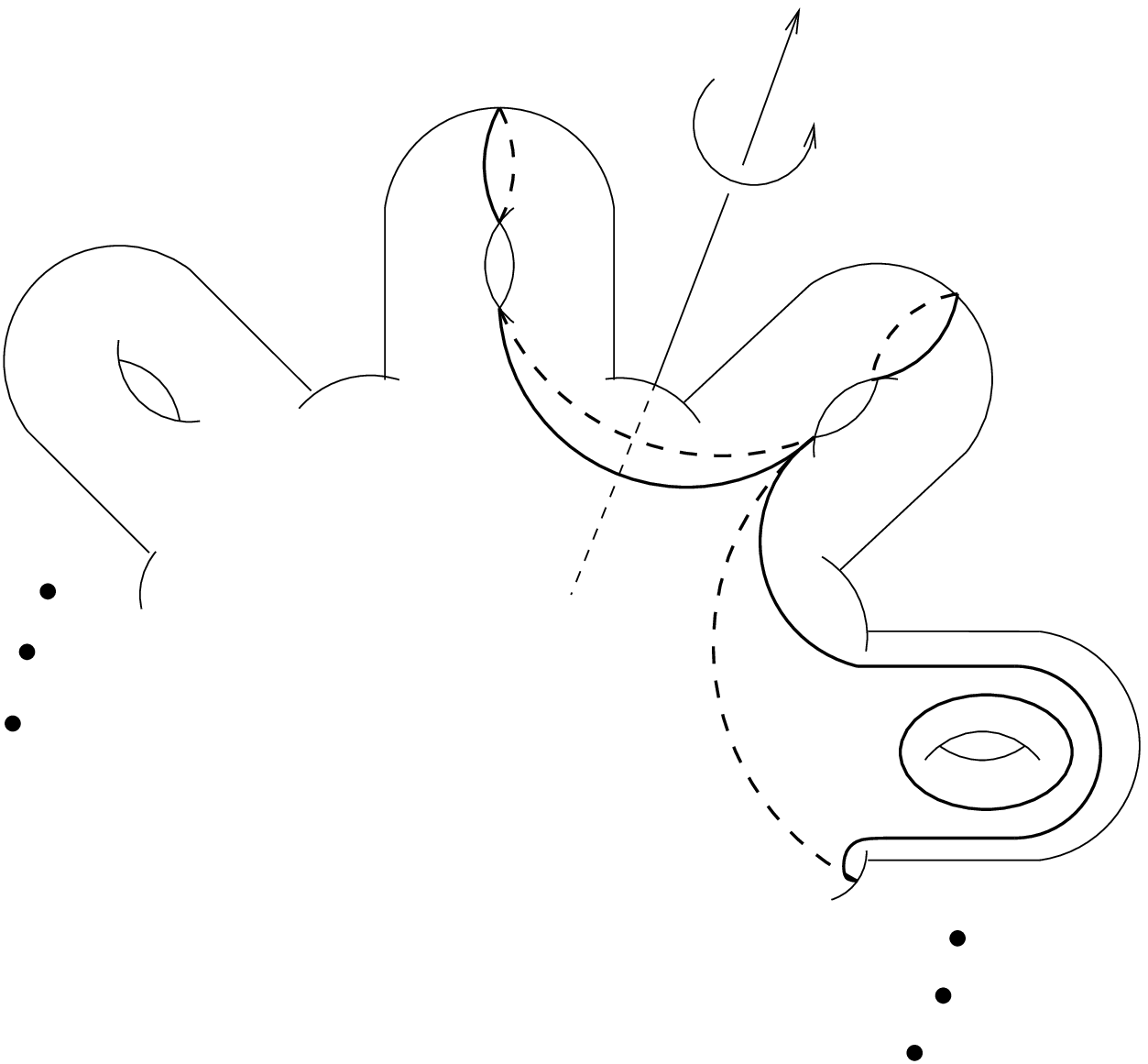}
$$
\caption{A good lantern embedded with $\alpha, \beta$ and $\gamma$
as three boundary curves and $x_1$ an interior curve.}
\label{good}
\end{figure}

We now show how to reduce this number.  Since, as previously noted, we
can replace $T_\beta$ with $T_{R_g^2 (\beta)}$ in our generating set,
we shall abuse notation by using $\beta$ to refer to the curve
${R_g^{2} (\beta)}$ for the remainder of the proof.  As shown in
Figure~\ref{good}, we can embed a good lantern in our surface which
contains $\alpha, \beta$, and $\gamma$ as three of its four boundary
components.  In order to motivate our choices of notation, we first
observe that, with this particular choice of a good lantern, we can
choose $x_1$ to be the curve labelled as such in Figure~\ref{good}.
If we then let $\alpha$ play the role of $a_1$, we see that $\rho_1
(a_1) = x_1$, and thus $\rho_1$ plays the role of $I_1$ in the proof
of Lemma~\ref{LH}.  

Now, in order to make further use of our work in
proving Lemma~\ref{LH}, we assign to the curve $\gamma$ the role of
$a_4$, so that $a_1$ and $a_4$ are separated by $x_1$ from $a_2$ and
$a_3$.  Of the two remaining curves, we now assign to $\beta$ the role
of $a_3$, and the remaining curve we label $a_2$.  We also choose
interior curves $x_2$ and $x_3$ so the labels of the good lantern of
Figure~\ref{good} match those of the pair swaps defined in the
previous section.  Thus we have
\begin{eqnarray}
T_\gamma = A_4 &=& (X_1 T_{\alpha}^{-1}) (X_2 A_2^{-1}) (X_3
T_\beta^{-1}) \label{eq:l1}
\\
&=&(X_1A_1^{-1}) (X_2 A_2^{-1}) (X_3 A_3^{-1}) \label{eq:l2} \\
&=&[ (X_1 I_1 X_1^{-1}) I_1 ] [ J_1 (X_1 I_1 X_1^{-1}) I_1 J_1 ]
[J_2(X_1 I_1 X_1^{-1}) I_1 J_2 ] \nonumber \\
&=& [ (X_1 \rho_1 X_1^{-1}) \rho_1 ] [ J_1 (X_1 \rho_1 X_1^{-1})
\rho_1 J_1 ] [J_2(X_1 \rho_1 X_1^{-1}) \rho_1 J_2 ] \nonumber
\end{eqnarray}

Given the two involutions $\rho_1,\rho_2$ used to generate $R_g$,
we have used only three new involutions to write $T_\gamma$,
namely $X_1 \rho_1 X_1^{-1}, J_1$, and $J_2$.  Now (\ref{eq:l1})
and ~(\ref{eq:l2}) can be rewritten, respectively, as
\begin{eqnarray}
T_\beta = A_3 &=& (X_1 T_{\alpha}^{-1}) (X_2 A_2^{-1}) (X_3
T_\gamma^{-1})
\nonumber\\
&=& (X_1A_1^{-1}) (X_2 A_2^{-1}) (X_3 A_4^{-1}) \label{eq:l4}
\end{eqnarray}

\begin{figure}
$$
\setlength{\unitlength}{0.05in}
\begin{picture}(0,0)(0,11)
\put(7,22.5){ {\bf $x_1$} }

\put(18.5,23){ {\bf $x_3$} }

\put(29,29){ {\bf $a_1$} }

\put(29,23.5){ {\bf $a_4$} }

\put(29,18){ {\bf $a_3$} }

\put(29,12){ {\bf $a_2$} }

\put(55,18){ {\bf $J_3$} }
\end{picture}
\includegraphics[width=3in]{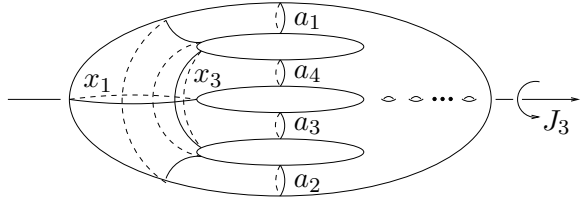}
$$
\caption{A pair swap $J_3$ which takes $(a_3, x_3)$ to $(a_4, x_3)$;
both $x_1$ and $x_3$ are interior lantern curves fixed by $J_3$.}
\label{twist2}
\end{figure}

We now observe that we can swap $(a_i,x_j)$ with $(a_k, x_j)$ for
any $i,j,k$ by using a slight variation on the embedding of
Figure~\ref{twist}. The involution $J_3$ shown in
Figure~\ref{twist2} interchanges $(a_3,x_3)$ with $(a_4,x_3)$,
as well as $(a_1,x_3)$ with $(a_2,x_3)$, $(a_1,x_1)$ with
$(a_2,x_1)$, and finally $(a_3,x_1)$ with $(a_4,x_1)$. Note that
in the first two cases, the $a_i$-curves that are interchanged
are not separated in the lantern by the fixed curve $x_3$, but in
the last two cases, the $a_i$-curves being swapped are separated by
the fixed curve $x_1$. Thus any possible configuration of such
pair swaps can be achieved, although a different embedding may be
required in order to align the $x_j$-curve to be fixed in the swap
either ``horizontally'' or ``vertically'', depending, respectively, on
whether the two $a_i$-curves one wishes to swap are separated by the
$x_j$-curve, or lie on the same side.

We now have that the product $J_3 J_2$ takes the pair $(x_1, a_1)$
to the pair $(x_3, a_4)$ and we can rewrite (\ref{eq:l4}) as
$$T_\beta = A_3 = (X_1A_1^{-1}) [J_1 (X_1 A_1)^{-1})
J_1] [(J_3 J_2) (X_1 A_1^{-1})(J_2 J_3)]$$ Since $(X_1A_1^{-1}) =
(X_1 \rho_1 X_1^{-1}) \rho_1$ as above, we see that we need just
one new involution, $J_3$, in order to generate $T_\beta$.
Similarly, we can find a pair-swap involution $J_4$ taking the
pair $(x_1, a_1)$ to $(x_1, a_4)$, so that we can write
\begin{eqnarray*}
T_\alpha = A_1 &=& (X_1 T_{\beta}^{-1}) (X_2 A_2^{-1}) (X_3 T_\gamma^{-1})\\
&=&(X_1 A_4^{-1}) (X_2 A_2^{-1}) (X_3 A_3^{-1})\\
&=&[J_4 (X_1 A_1^{-1}) J_4] [J_1 (X_1 A_1^{-1}) J_1] [J_2 (X_1 A_1^{-1})J_2 ]
\end{eqnarray*}
\noindent with $J_4$ being the only new involution required to
write $T_\gamma$. 

Thus our count for the number of involutions used to generate each of
the four elements $R_g, T_\gamma, T_\beta, T_\alpha$ stands,
respectively, at $2 + 3 + 1 + 1 = 7$ involutions.  Finally, as pointed
out to us by M. Kassabov, the involution $J_4$ is in fact redundant,
and can be taken to be the product $J_2 J_3 J_2$.  Thus we can generate 
with $6$ involutions, as claimed.  \endproof

\section{Final remarks}
\label{section:final}

Recall that an {\em Artin group} $A$ is generated by elements $x_1,
\ldots, x_n$, subject to the following relations.  Let $m_{s,t} =
\{2,3,4, \ldots, \infty \}$.  Then for each $s \neq t$, the elements
$x_s$ and $x_t$ satisfy the relations:

$$
\begin{array}{ll} 
(x_s x_t)^{m_{s,t}/2} = (x_t
x_s)^{m_{s,t}/2} &\mbox{if $m_{s,t}$ is even}\\
&\\
(x_sx_t)^{(m_{s,t}-1)/2} = (x_t x_s)^{(m_{s,t}-1)/2}&
\mbox{if $m_{s,t}$ is odd}
\end{array}
$$

Here $m_{s,t} = \infty$ means that there is no relation between
the generators $x_s$ and $x_t$.  The {\em Coxeter group} $\bar{A}$
associated to an Artin group $A$ is the quotient of $A$ by the extra
relation $x_i^2 = 1$ for all $i = 1, \ldots n$. We thank Nathan
Broaddus for pointing out the following corollary of our main theorem.

\begin{corollary}[\boldmath$\mgb$ is a Coxeter quotient]
\label{coxeter} 
For $g \geq 3, b = 0$, and
for $g \geq 4, b = 1$, the group $\mgb$ can be realized as a
quotient of a Coxeter group on 6 generators.
\end{corollary}

There has been some effort to understand various mappings of Artin
groups both into and onto mapping class groups.  For example,
Wajnryb studied classes of Artin groups which do not inject into
mapping class groups in \cite{Wa2}, while Matsumo \cite{Mat} as
well as Labruere and Paris \cite{LP} have given explicit
presentations of mapping class groups as quotients of certain
Artin groups. However, we are not aware of similar investigations
of the relationship between Coxeter groups and mapping class
groups.

The results of this paper beg several further questions.  For
example, besides the previously raised question of whether we can
do better than $6$ involutions, we can also
ask whether there exists a constant $C$, perhaps such that $C =
C(g)$, so that every element of $\mgb$ can be written as a product
of at most $C$ torsion elements \footnote{After initial
circulation of this paper, M. Korkmaz \cite{Ko2} and D. Kotschick \cite{Kot} 
answered the latter question in the negative by proving that no such $C$
exists.  Both arguments build on a theorem of
Endo-Kotschick \cite{EK}.} . Furthermore, what kinds of relations
exist amongst the torsion (or involution) generators?  In
particular, what kind of Coxeter groups arise in the context of
Corollary~\ref{coxeter}?

\noindent
Tara E. Brendle:\\
Dept. of Mathematics, Cornell University\\
310 Malott Hall\\
Ithaca, NY  14853\\
E-mail: brendle@math.cornell.edu
\medskip

\noindent
Benson Farb:\\
Dept. of Mathematics, University of Chicago\\
5734 University Ave.\\
Chicago, Il 60637\\
E-mail: farb@math.uchicago.edu

\end{document}